\documentclass [12pt] {amsart}
\usepackage{times}
\usepackage{latexsym,amssymb, amsmath}
\usepackage[margin=2.5cm]{geometry}

\def\pf{\noindent\emph{Proof: }}
\def\stop{\hfill$\Box$}

\usepackage{amsmath, amsfonts, amssymb, amsthm}
\newtheorem{thm}{Theorem}

\newtheorem{lemma}[thm]{Lemma}

\newtheorem{example}{Example}

\newtheorem{conj}{Conjecture}
\numberwithin{thm}{section}

\begin{document}

\title{Uniqueness Theorems for Ordinary Differential Equations with H\"older Continuity}

\author{Yifei Pan}
\address{Department of Mathematical Sciences\\
 Indiana University-Purdue University Fort Wayne\\
 Fort Wayne, Indiana 46805}
\email{pan@ipfw.edu}

\author {Mei Wang}
\address {Department of Statistics\\University of Chicago \\Chicago, Illinois 60637}
\email{meiwang@galton.uchicago.edu}

\author { Yu Yan  }
\address {Department of Mathematics and Computer Science\\ Huntington University\\Huntington, Indiana 46750}
\email {yyan@huntington.edu}

\begin{abstract}
We study ordinary differential equations of the type $u^{(n)}(t)=f(u(t))$ with initial conditions $u(0) =  u'(0)  =  \cdots   = u^{(m-1)}(0) =  0  $ and $u^{(m)}(0) \neq 0$ where $m  \geq n$, no additional assumption is made on $f$.   We establish some uniqueness results and show that $f$ is always H\"older continuous.
\end{abstract}

\maketitle
\newtheorem{Thm}{Theorem}
\newtheorem{Lemm}{Lemma}
\newtheorem{Cor}{Corollary}

\section {Introduction}


The question of finding criteria for the uniqueness of solutions has been a constant theme in the study of ordinary differential equations for a very long time, and a wealth of results have been established.  The most quoted one in textbooks is perhaps the Lipschitz uniqueness theorem, which states that in the equation $y^{(n)}(x)=f(x, y, y',..., y^{(n-1)})$, if the function $f(x,z_1,z_2,...z_n)$ is Lipschitz continuous with respect to $z_1,z_2,...z_n$, then the initial value problem has a unique local solution.  Generally speaking, to ensure the uniqueness of solutions to an ODE, we need to assume some condition on the function $f$ besides continuity, the Lipschitz condition being one example.  Most of the research in this topic has been devoted to finding the appropriate condition, and there are many nice results such as the classical theorems by Peano, Osgood, Montel-Tonelli, and Nagumo.  The monograph \cite{AAL_ODE} provides an extensive and systematic treatment of the available results.

\vspace{.1in}

In this paper, we approach the uniqueness problem from a different perspective and relate it to the unique continuation problem.  We study autonomous ODEs of the type  $u^{(n)}(t)=f \left (u \left (t \right ) \right )$ where $u \in C^{\infty}([0,1])$ and no additional assumption is made on the function $f$.

\vspace{.1in}
If we assume the initial conditions $u(0) =  u'(0)  =  \cdots   = u^{(n-1)}(0) = \,\, 0,$ the solution is not unique.  The following is a trivial example.

\begin{example}
$u(t)=t^3$ satisfies $u''(t)=6u^{\frac{1}{3}}$ and $u(0)=u'(0)=0$.  Another solution to this initial value problem is $u \equiv 0$.
\end{example}

\noindent
It is no surprise that uniqueness fails in this example because the function $f(u)=6u^{\frac{1}{3}}$ has fairly strong singularity at $0$.  From another perspective, this example shows that if a solution and its derivatives up to order $n-1$ all vanish at 0, it is not guaranteed to be the zero function.  On the other hand, however, even if all its derivative vanish at 0, the solution still may not be identically 0.

\begin{example}
The function
\begin{displaymath}
u(t) = \left\{
\begin{array}{lr}
e^{-\frac{1}{t}} &  \hspace{.1in}  0 < t\leq 1\\
0 &  \hspace{.1in} t=0
\end{array}   \right.
\end{displaymath}

\noindent
is in $C^{\infty} \left ( [0,1] \right )  $ and
\begin{equation}
\label{eq:infinite_example}
u^{(k)}(0)=0  \hspace{.3in} \text{ for all} \hspace{.2in} k \in \mathbf{N}.
\end{equation}
Let

\begin{displaymath}
f(s) = \left\{
\begin{array}{lr}
s(\ln s)^2 &  \hspace{.1in}  s>0\\
0 &  \hspace{.2in} s=0
\end{array}   \right.
\end{displaymath}

\noindent
Then $u(t)$ satisfies the equation $$u'=f(u).$$   However, this equation has another solution $u \equiv 0$, which also satisfies (\ref{eq:infinite_example}).

\end{example}

\noindent
This function $u(t)$ is also a classical example in the study of the unique continuation problem, that is, when can we conclude that a function is locally identically zero if its derivatives all vanish at a point.  One result along this line is given in \cite{Pan-Wang_1}.

\vspace{.05in}

\begin{thm}
\label{thm:Pan-Wang_flatness}
(\cite{Pan-Wang_1}) Let $g(x) \in C^{\infty}([a,b])$, $ 0 \in [a,b]$, and
\begin{equation}
\label{eq: 2nd_derivative_condition}
|g^{(n)}(x)| \leq C \sum_{k=0}^{n-1} \frac{\left | g^{(k)}(x) \right |}{|x|^{n-k}}, \hspace{.5in} x \in [a,b]
\end{equation}
for some constant $C$ and some $n \geq 1$. Then
$$g^{(k)}(0)=0 \hspace{.2in} \forall k \geq 0$$ implies $$g \equiv 0 \hspace{.1in} \text{on} \,\, [a,b]$$
\end{thm}

\vspace{.1in}

\noindent
The order of singularity of $|x|$ at 0 in (\ref{eq: 2nd_derivative_condition}), i.e. $n-k$, is sharp, as one example in \cite{Pan-Wang_1} shows.  This theorem is crucial to the proof of our main theorem below.

\vspace{.1in}
The previous two examples suggest that to guarantee uniqueness near 0, the solution needs to vanish to sufficiently high order, but not to the infinite order.  So we assume that it satisfies the initial conditions $u(0) =  u'(0)  =  \cdots   = u^{(m-1)}(0) = \,\, 0 $ and $ u^{(m)}(0) = a \neq 0$, where $m\geq n$.  That is, the order of the lowest non-vanishing derivative of $u$ at $0$ is no less than the order of the equation.  From the equation it is not difficult to see that $f$ is differentiable away from $0$, however, it is not differentiable at $0$, as shown by Example 1 with $m=3$.

\vspace{.1in}

\noindent
Due to the lack of information about the regularity of $f$, the available uniqueness theory no longer applies to this type of equations.  We will show that because $u$ has sufficiently high vanishing order at $0$, such solutions are unique near $0$.   Specifically, we have the following result.

\begin{thm}
\label{thm:uniqueness-ODE}
Let $u(t) \in C^{\infty}([0,1))$ be a solution of the differential equation
\begin{equation}
\label{eq:main}
u^{(n)}(t)=f \left (u \left (t \right ) \right ),
\end{equation}
 where $n \geq 1$ and  $f$ is a function.  Assume that $u$ satisfies
\begin{equation}
\label{eq:initial}
u(0) \,\, = \,\, u'(0) \,\, = \,\, \cdots  \,\, = u^{(m-1)}(0) \,\, = \,\, 0 \hspace{.1in} \text{and} \hspace{.1in} u^{(m)}(0) = a \neq 0
\end{equation}
 with $m \geq n$. Then such a solution $u(t)$ is unique for $t$ near 0.

\end{thm}

\vspace{.05in}

\vspace{.05in}
\noindent
The proof of Theorem \ref{thm:uniqueness-ODE} is carried out in two steps.  First, we show the following result concerning the derivatives of $u$ at $0$.

\vspace{.05in}

\begin{lemma}
\label{lem:step1_Taylor_coefficient}
Let $u(t)$ be a solution that satisfies Equations (\ref{eq:main} ) and  (\ref{eq:initial}). The derivative of $u$ at 0 of any order equal to or higher than $m$, that is, $u^{(k)}(0) $ for any   $ k \geq m$, depends only on $m,n$, and the behavior of the function $f$ near 0.
\end{lemma}

\vspace{.05in}

\noindent
In the second step, suppose there are two solutions $u$ and $v$, both satisfying (\ref{eq:main}) and  (\ref{eq:initial}), then by the previous lemma the function $u(t)-v(t)$ and all its derivatives vanish at $0$.  Making use of Theorem \ref{thm:Pan-Wang_flatness}, we can show that $u-v \equiv 0$.

\vspace{.05in}

\noindent
Typically, for an $n$-th order ODE we need only $n$ initial conditions.  This theorem shows that in some sense, the lack of information about $f$ can be compensated by assuming additional derivative information at the initial point.

\vspace{.1in}

\vspace{.1in}

Interestingly, it turns out that the solution is unique as long as the vanishing order is no less than the order of the equation, but the actual vanishing order and the value of the lowest non-zero derivative are not essential.

\vspace{.05in}

\begin{thm}
\label{thm:u1_u2}
Suppose $u_1$ and $u_2$ are two solutions of Equation (\ref {eq:main}) and they satisfy
\begin{equation}
\label{eq:initial_u1}
u_1(0) \,\, = \,\, u_1'(0) \,\, = \,\, \cdots  \,\, = u_1^{(m-1)}(0) \,\, = \,\, 0, \hspace{.2in}  u_1^{(m)}(0) = a \neq 0
\end{equation}

\noindent
and
\begin{equation}
\label{eq:initial_u2}
u_2(0) \,\, = \,\, u_2'(0) \,\, = \,\, \cdots  \,\, = u_2^{(l-1)}(0) \,\, = \,\, 0, \hspace{.2in}  u_2^{(l)}(0) = b \neq 0
\end{equation}
where $m,l \geq n$.  Then $m=l$, $a=b$ and $u_1 \equiv u_2$ for small $t$.
\end{thm}

\vspace{.05in}

The proof of Lemma \ref{lem:step1_Taylor_coefficient} will be given in Section \ref{sec:step1_proof}, and the proof of Theorems \ref{thm:uniqueness-ODE} and \ref{thm:u1_u2} will be given in Section \ref{sec:main_pf}.

\vspace{.05in}

Naturally, we would like to ask if the result still holds if one of the solutions in Theorem \ref{thm:u1_u2} has vanishing order lower than $n$, the order of the equation.  We propose the following conjecture.

\vspace{.05in}
\begin{conj}
Suppose Equation (\ref{eq:main}) has a solution $u(t)$ which satisfies (\ref{eq:initial}) with $m \geq n$. Then it cannot possess another solution $v(t)$ which satisfies initial conditions
$$
v(0) \,\, = \,\, v'(0) \,\, = \,\, \cdots  \,\, = v^{(l-1)}(0) \,\, = \,\, 0 \hspace{.1in} \text{and} \hspace{.1in} v^{(l)}(0) = b \neq 0,
$$ where $l<n$.
\end{conj}

\vspace{.05in}
\noindent
We can show that this conjecture is true if $m=n+1$ or $l$ and $m$ are relatively prime.  However, there are some difficulties in the general case and we have not been able to prove the full conjecture.

\vspace{.1in}
\noindent

 Although in Theorems \ref{thm:uniqueness-ODE} and \ref{thm:u1_u2} we do not need to make any assumptions about the function $f$, we can actually obtain interesting information about it.  Suppose there is a function $u(t)$ which satisfies Condition (\ref {eq:initial}), then as shown in Section \ref{sec:step1_proof}, locally $t$ can be expressed as a function of $u$, therefore we can express $u^{(n)}(t)$ locally as a function $f$ of $u$, so $u^{(n)}(t)=f(u)$.  The next theorem shows that the function $f$ is H\"older continuous in an interval $[0,\delta]$ for small $\delta >0$.

\begin{thm}
\label{thm:f_ho"lder}
Suppose a function $u(t)$ satisfies Condition (\ref {eq:initial}), then Equation (\ref {eq:main}) holds for some function $f$ where $m \geq n$, and there is a constant $\delta >0$ such that $f$ is uniformly H\"older continuous in the interval $[0,\delta]$.
\end{thm}

\noindent
This theorem is proved after Theorem \ref{thm:u1_u2} in Section \ref{sec:main_pf}.

\vspace{.1in}
A summary of Theorems \ref{thm:uniqueness-ODE} and \ref{thm:f_ho"lder} is that any smooth function of finite order vanishing at $0$ is a unique local solution of a differential equation in the form of (\ref {eq:main}), where $f$ is differentiable in the interior and uniformly H\"older continuous up to the boundary.

\vspace{.1in}
In Theorems \ref{thm:uniqueness-ODE} and \ref{thm:u1_u2}, the high order vanishing condition (\ref{eq:initial}) allows us to obtain uniqueness results without any extra assumption on $f$.  This phenomenon is only found in autonomous equations like (\ref {eq:main}).  For general ODEs of the form $\frac{d^nu}{dt^n}=f \left (t,u,\frac{du}{dt},...,\frac{d^{n-1}u}{dt^{n-1}}  \right )$, such results cannot be expected because there is more than one expression for $f$.  For example,

\begin{example}
\label{example_u'}
Let $u(t)=t^4$.  It satisfies the initial conditions $$u(0)=u'(0)=u''(0)=u^{(3)}(0)=0 \hspace{.2in}  \text{and} \hspace{.2in} u^{(4)}(0)=24.$$  Its derivatives are $$u'(t)=4t^3, \hspace{.2in}  \text{and} \hspace{.2in} u''(t)=12t^2.$$  We can express $u''$ as a function $f$ of $u$ and $u'$ in different ways such as
$$u''(t)=\frac{192u^2}{(u')^2},$$  or $$u''(t)=\frac{3(u')^2}{4u},$$ or $$u''(t)=12 \left ( \frac{1}{4}uu' \right )^{\frac{2}{7}},$$ or $$u''(t)=12 u^{\frac{1}{4}} \left ( \frac{1}{4}u'\right )^{\frac{1}{3}}.$$
\end{example}
\vspace{.1in}

\noindent
In the first two equations $f$ is not continuous at the origin, while in the last two equations $f$ is H\"older continuous at the origin.  This simple example shows that the function $f$ can be expressed in various ways and that in order to study the uniqueness we need to impose very specific assumptions on $f$.

\vspace{.1in}
Our work was motivated by \cite{Li-Ni_1} in which Li and Nirenberg studied a similar second order PDE:  $\Delta u=f(u)$, where $u=u(t,x)\in C^{\infty}(\mathbf{R}^{k+1})$ has a non-vanishing partial derivative at 0 which can be expressed in the form $u(t,x)=at^m+O(t^{m+1})$, $a \neq 0$, $t \in \mathbf{R}$ and $ x \in \mathbf{R}^k $.  They showed that if two solutions $u$ and $v$ satisfy $u\geq v$, then $u \equiv v$.  Theorem \ref{thm:uniqueness-ODE} can be viewed as an improvement of their result in the one-dimensional case to arbitrary order and without the comparison condition $u\geq v$.

\vspace{.2in}

\section {The Proof of Lemma \ref{lem:step1_Taylor_coefficient}  }
\label{sec:step1_proof}

\noindent
Without loss of generality, we can assume that $a>0$.
\vspace{.1in}

\noindent

\begin{itemize}
	\item First, we show that $a$, the $m$-th derivative of $u$ at 0, only depends on $m,n$, and the function $f$.
\end{itemize}

\vspace{.1in}

\noindent
Define
\begin{equation}
\label{eq:tilde-x-defn}
\tilde{x}= \left (  \frac{u}{a}\right )^{\frac{1}{m}}.
\end{equation}

\noindent
 Then $$\tilde{x} \,\, = \,\, \left (  t^m+O \left (t^{m+1} \right ) \right )^{\frac{1}{m}} \,\, = \,\, t \left ( 1+O(t) \right ).$$
 This implies that

\begin{equation}
\label{eq:x-tilde-over-t}
\frac{\tilde{x}}{t}  \rightarrow 1 \hspace{.2in} \text{as}  \hspace{.2in} t  \rightarrow 0 .
\end{equation}

\noindent
We can also write $u=a\tilde{x}^m$.  Taking the derivative with respect to $t$, we get

\begin{eqnarray*}	
	\frac{du}{dt} & = & am\tilde{x}^{m-1} \frac{d\tilde{x}}{dt}\\
	amt^{m-1}+O(t^m) & = &  am\tilde{x}^{m-1} \frac{d\tilde{x}}{dt}\\
	\frac{t^{m-1}}{\tilde{x}^{m-1}}+O \left (  \frac{t^m}{\tilde{x}^{m-1}}\right ) & = &\frac{d\tilde{x}}{dt}
\end{eqnarray*}

\vspace{.1in}
\noindent
In the second equation above and in the analysis that follows, we formally differentiate the Taylor expansion with the big-$O$ notation.  A detailed discussion of this differentiation is provided in the Appendix.

\vspace{.1in}
\noindent
In light of (\ref{eq:x-tilde-over-t}), it follows that

\begin{equation}
\label{eq:dx-tilde-dt}
\frac{d\tilde{x}}{dt} | _{t=0}=1.
\end{equation}

\noindent
By the Inverse Function Theorem, $t$ can be expressed as a function of $\tilde{x}$: $t=\tilde{x}+O\left (  \tilde{x}^2 \right ) .$  Then

$$ t^{m-n} \,\, = \,\,  \left (  \tilde{x}+O\left (  \tilde{x}^2 \right )  \right )^{m-n} \,\, = \,\, \tilde{x} ^{m-n} \left ( 1+O \left ( \tilde{x}\right ) \right ) .$$  Similarly, $ t^{m-n+1} \,\, = \,\, \tilde{x} ^{m-n+1} \left ( 1+O \left ( \tilde{x}\right ) \right ). $ Thus

\begin{eqnarray}
\label{eq:for-use-in-main-thm}
f(u) & = & u^{(n)} \nonumber \\
& = & am(m-1)\cdots(m-n+1)t^{m-n}+O \left ( t^{m-n+1}\right )  \nonumber \\
& = & am(m-1)\cdots(m-n+1)\tilde{x}^{m-n}+O \left ( \tilde{x}^{m-n+1}\right ) \nonumber \\
& = & am(m-1) \cdots(m-n+1)\left ( \frac{u}{a} \right ) ^ {\frac{m-n}{m}}+ O \left ( u^\frac{m-n+1}{m}\right )  \hspace{.3in} \text {by Equation (\ref {eq:tilde-x-defn})} \nonumber \\
& = & a^\frac{n}{m} m(m-1)\cdots(m-n+1) u^{\frac{m-n}{m}} +  O \left ( u^\frac{m-n+1}{m}\right )
\end{eqnarray}

\vspace{.1in}
\noindent
Therefore,

\begin{equation}
\label{eq:a-determined}
a^{\frac{n}{m}} = \lim _{u \to 0} \frac{f(u)}{m(m-1) \cdots(m-n+1)u^{\frac{m-n}{m}}}.
\end{equation}
This shows that $a$ is completely determined by $m,n$, and the behavior of the function $f$ near 0.

\vspace{.2in}

\begin{itemize}
	\item Next, we show that the $(m+1)$-th derivative of $u$ at 0 also only depends on $m,n$, and $f$.
\end{itemize}

\noindent
Write $u$ as

\begin{equation}
\label{eq:u_m+1}
u(t)=at^m+a_{m+1}t^{m+1}+O(t^{m+2}).
\end{equation}
We will show that $a_{m+1}$ only depends on $m,n$, and the behavior of $f$ at $0$.

\noindent
Express $t$ as

\begin{equation}
\label{eq:t-2nd-order}
t=\tilde{x}+b_2\tilde{x}^2+O \left (  \tilde{x}^3\right ).
\end{equation}

\noindent
We would like to obtain an expression for $b_2$ in terms of the derivatives of $u$ at $0$.  To do this, we take the derivative with respect to $t$ on both sides of $\frac{d\tilde{x}}{dt} \cdot \frac{dt}{d\tilde{x} } = 1$.  By the Product Rule and Chain Rule, we have

\begin{eqnarray}
\label{eq:x-tilde-t-1st}
 \frac{d}{dt} \left ( \frac{d\tilde{x}}{dt} \right ) \cdot \frac{dt}{d\tilde{x}} + \frac{d\tilde{x}}{dt} \cdot \frac{d}{dt} \left ( \frac{dt}{d\tilde{x}} \right ) & = & 0 \nonumber \\
\frac{d^2\tilde{x}}{dt^2}\cdot \frac{dt}{d\tilde{x}} + \frac{d\tilde{x}}{dt} \cdot \left (\frac{d^2t}{d\tilde{x}^2} \cdot \frac{d\tilde{x}}{dt} \right ) & = & 0 \nonumber \\
\frac{d^2\tilde{x}}{dt^2}\cdot \frac{dt}{d\tilde{x}} + \left ( \frac{d\tilde{x}}{dt} \right )^2 \cdot  \frac{d^2t}{d\tilde{x}^2}   & = & 0
\end{eqnarray}

\noindent
We would like to evaluate Equation (\ref{eq:x-tilde-t-1st}) at $t=0$.

\noindent
From
\begin{eqnarray*}
	\tilde{x} & = & \left (  \frac{u}{a}\right )^{\frac{1}{m}}\\
	& = &  \left ( t^m+ \frac{a_{m+1}}{a} t^{m+1}+O \left (t^{m+2} \right ) \right ) ^{\frac{1}{m}}\\
	& = &  t \left [ 1+ \frac{a_{m+1}}{a} t + O (t^2) \right ] ^{\frac{1}{m}}\\
	& = &  t \left [ 1+ \frac{1}{m} \left (\frac{a_{m+1}}{a} t + O (t^2) \right )  + \frac{1}{2} \cdot\frac{1}{m} \left ( \frac{1}{m}-1 \right)  \left (\frac{a_{m+1}}{a} t + O (t^2) \right )^2 +O(t^3)\right ]\\
	& = & t+\frac{a_{m+1}}{ma} t^2+O(t^3),
\end{eqnarray*}
we know that

\begin{equation}
\label{eq:d-tilde-x-dt=2nd}
\frac{d^2\tilde{x}}{dt^2} | _{t=0}=\frac{2a_{m+1}}{ma}.
\end{equation}

\noindent
From (\ref{eq:t-2nd-order}) we know that $$\frac{dt}{d\tilde{x}} |_{t=0}=1 \hspace{.3in} \text { and } \hspace{.3in} \frac{d^2t}{d\tilde{x}^2} |_{t=0}=2b_2.$$  Thus if we evaluate Equation (\ref{eq:x-tilde-t-1st}) at $t=0$, we get $$ \frac{2a_{m+1}}{ma} \cdot 1 + 1 \cdot 2b_2=0,$$ therefore

\begin{equation}
\label{eq:b2}
b_2=-\frac{a_{m+1}}{ma}.
\end{equation}

\noindent
Now, from Equation (\ref{eq:t-2nd-order}) we have

\begin{eqnarray*}
t^{m-n} & = & \left (  \tilde{x}+b_2\tilde{x}^2+O \left (  \tilde{x}^3\right ) \right ) ^{m-n} \\
& = & \tilde{x}^{m-n} \left [  1 + b_2\tilde{x}+O \left (  \tilde{x}^2\right ) \right ] ^{m-n} \\
& = & \tilde{x}^{m-n} \left [  1 + (m-n) \left ( b_2\tilde{x}+O \left (  \tilde{x}^2 \right ) \right ) + O \left (  \tilde{x}^2 \right ) \right ]  \\
& = & \tilde{x}^{m-n} + (m-n)b_2 \tilde{x}^{m-n+1} + O \left (  \tilde{x}^{m-n+2} \right ).
\end{eqnarray*}
Similarly
\begin{eqnarray*}
t^{m-n+1} & = &  \tilde{x}^{m-n+1} + (m-n+1)b_2 \tilde{x}^{m-n+2} + O \left (  \tilde{x}^{m-n+3} \right ).\\
t^{m-n+2} & = & O \left (  \tilde{x}^{m-n+2} \right ).
\end{eqnarray*}

\noindent
Then from Equation (\ref{eq:u_m+1}) and the above expressions for the powers of $t$, we have

\begin{eqnarray*}
u^{(n)} & = & am(m-1)\cdots(m-n+1) t^{m-n} \\
& & +  a_{m+1}(m+1)m \cdots(m-n+2) t^{m-n+1} + O \left ( t^{m-n+2} \right ) \\
& = & am(m-1) \cdots(m-n+1)\left [  \tilde{x}^{m-n} + (m-n)b_2 \tilde{x}^{m-n+1} + O \left (  \tilde{x}^{m-n+2} \right ) \right ] \\
& & + a_{m+1}(m+1)m \cdots(m-n+2) \left [ \tilde{x}^{m-n+1} + (m-n+1)b_2 \tilde{x}^{m-n+2} + O \left (  \tilde{x}^{m-n+3} \right ) \right ] \\
& & + O \left (\tilde{x}^{m-n+2} \right ) \\
& = & am(m-1)  \cdots(m-n+1)\tilde{x}^{m-n}  \\
& & + \left [ am(m-1)\cdots (m-n)b_2 + a_{m+1}(m+1)m \cdots(m-n+2)  \right ] \tilde{x}^{m-n+1} + O \left (\tilde{x}^{m-n+2} \right )
\end{eqnarray*}

\noindent
Thus by $f(u)=u^{(n)}$ and (\ref {eq:tilde-x-defn} ), we have

\begin{eqnarray}
f(u) & = & am(m-1)  \cdots(m-n+1)\tilde{x}^{m-n}  \nonumber \\
& & + \Big [ am(m-1)\cdots (m-n)b_2 + a_{m+1}(m+1)m \cdots(m-n+2)  \Big ] \tilde{x}^{m-n+1} + O \left (\tilde{x}^{m-n+2} \right )  \nonumber \\
& = & am(m-1)  \cdots(m-n+1) \left ( \frac{u}{a} \right ) ^{\frac{m-n}{m}}  \nonumber \\
& & + \Big [ am(m-1)\cdots (m-n)b_2 + a_{m+1}(m+1)m \cdots(m-n+2)  \Big ] \left ( \frac{u}{a} \right ) ^{\frac{m-n+1}{m}} \nonumber \\
& & + O \left (\left ( \frac{u}{a} \right ) ^{\frac{m-n+2}{m}} \right )  \nonumber \\
& = & a^{\frac{n}{m}}m(m-1) \cdots(m-n+1) u ^{\frac{m-n}{m}} \nonumber \\
& & + \Big [ a^{\frac{n-1}{m}}m(m-1)\cdots (m-n)b_2 + a^{\frac{n-m-1}{m}}a_{m+1}(m+1)m \cdots(m-n+2) \Big ] u^{\frac{m-n+1}{m}}  \nonumber \\
& & + O \left ( u  ^{\frac{m-n+2}{m}}\right ) \nonumber
\end{eqnarray}

\noindent
This means that

\begin{eqnarray*}
& & a^{\frac{n-1}{m}}m(m-1)\cdots (m-n)b_2 + a^{\frac{n-m-1}{m}}a_{m+1}(m+1)m \cdots(m-n+2) \\
& = &\lim _{u \to 0} \frac{f(u)-a^{\frac{n}{m}}m(m-1) \cdots(m-n+1) u ^{\frac{m-n}{m}}}{u^{\frac{m-n+1}{m}}}.
\end{eqnarray*}

\noindent
By (\ref{eq:b2}), this can be written as
\begin{eqnarray*}
& & a^{\frac{n-1}{m}}m(m-1)\cdots (m-n) \left ( -\frac{a_{m+1}}{ma} \right ) + a^{\frac{n-m-1}{m}}a_{m+1}(m+1)m \cdots(m-n+2) \\
& = &\lim _{u \to 0} \frac{f(u)-a^{\frac{n}{m}}m(m-1) \cdots(m-n+1) u ^{\frac{m-n}{m}}}{u^{\frac{m-n+1}{m}}}.
\end{eqnarray*}

\noindent
After collecting similar terms, we get
\begin{eqnarray*}
& & a_{m+1}a^{\frac{n-m-1}{m}}\left [ (m+1)m\cdots(m-n+2) - (m-1)(m-2)\cdots (m-n) \right ] \\
& = &  \lim _{u \to 0} \frac{f(u)-a^{\frac{n}{m}}m(m-1) \cdots(m-n+1)  u ^{\frac{m-n}{m}}}{u^{\frac{m-n+1}{m}}}.
\end{eqnarray*}

\noindent
Consequently,
\begin{eqnarray}
\label{eq:a_m+1_determined}
 a_{m+1}  & = & \frac{a^{\frac{m-n+1}{m}} \cdot \left ( \displaystyle \lim _{u \to 0} \frac{f(u)-a^{\frac{n}{m}}m(m-1) \cdots(m-n+1) u ^{\frac{m-n}{m}}}{u^{\frac{m-n+1}{m}}} \right ) }{(m+1)m\cdots(m-n+2) - (m-1)(m-2)\cdots (m-n)}.
\end{eqnarray}

\vspace{.1in}

\noindent
Since we have proved that $a$ only depends on $m,n$, and $f$, Equation (\ref {eq:a_m+1_determined}) shows that $a_{m+1}$ is also completed determined by $m,n$, and the behavior of $f$ near 0.  By (\ref {eq:b2}), this also shows that $b_2$ depends only on $m,n$, and $f$.

\vspace{.1in}

\begin{itemize}
	\item Then, we will use mathematical induction to show that all the derivatives of $u$ at 0 of order higher than $m$ are completely determined by $m,n$, and $f$.
\end{itemize}

\vspace{.1in}

\noindent
Express $u$ and $t$ as
\begin{equation}
\label{eq:u_m+k}
u(t)=at^m + a_{m+1}t^{m+1} + \cdots +  a_{m+k}t^{m+k} +  a_{m+k+1}t^{m+k+1} + O \left ( t^{m+k+2}\right ),
\end{equation}

\begin{equation}
\label{eq:t_k+2}
t = \tilde{x} + b_2 \tilde{x}^2 + \cdots + b_{k+1} \tilde{x}^{k+1} + b_{k+2} \tilde{x}^{k+2} + O \left (  \tilde{x}^{k+3} \right )
\end{equation}

\vspace{.1in}
\noindent
Suppose that for $k \geq 1$, $a$, $a_{m+1}$, ..., $a_{m+k}$, $b_2$,..., $b_{k+1}$ are all determined only by by $m,n$, and $f$, we will show that $a_{m+k+1}$ and $b_{k+2}$ also are determined only by by $m,n$, and $f$.

\vspace{.1in}

\noindent
We start by obtaining an expression for $b_{k+2}$ in terms of $a_{m+1}$, ..., $a_{m+k}$ and $a_{m+k+1}$.

\vspace{.1in}

\noindent
Taking the derivative with respect to $t$ on both sides of (\ref {eq:x-tilde-t-1st}), we obtain

\begin{eqnarray}
\label {eq:x-tilde-t-2nd}
0 & = & \frac{d^3\tilde{x}}{dt^3}\cdot \frac{dt}{d\tilde{x}} + \frac{d^2\tilde{x}}{dt^2} \cdot \left ( \frac{d^2t}{d\tilde{x}^2} \cdot \frac{d\tilde{x}}{dt} \right ) + 2 \cdot \frac{d\tilde{x}}{dt}\cdot\frac{d^2\tilde{x}}{dt^2}\cdot\frac{d^2t}{d\tilde{x}^2} + \left ( \frac{d\tilde{x}}{dt} \right )^2 \cdot \left ( \frac{d^3t}{d\tilde{x}^3} \cdot \frac{d\tilde{x}}{dt} \right ) \nonumber \\
0 & = & \frac{d^3\tilde{x}}{dt^3}\cdot \frac{dt}{d\tilde{x}} + 3 \cdot \frac{d\tilde{x}}{dt}\cdot\frac{d^2\tilde{x}}{dt^2}\cdot\frac{d^2t}{d\tilde{x}^2} + \left ( \frac{d\tilde{x}}{dt} \right )^3\cdot \frac{d^3t}{d\tilde{x}^3}
\end{eqnarray}

\vspace{.1in}
\noindent
Taking the derivative of both sides of (\ref {eq:x-tilde-t-2nd}), we get

\begin{eqnarray}
\label {eq:x-tilde-t-3rd}
 0 & = & \left ( \frac{d^4\tilde{x}}{dt^4}\cdot \frac{dt}{d\tilde{x}} + \frac{d^3\tilde{x}}{dt^3} \cdot \frac{d^2t}{d\tilde{x}^2} \cdot \frac{d\tilde{x}}{dt} \right ) + 3 \Bigg ( \frac{d^2\tilde{x}}{dt^2}\cdot\frac{d^2\tilde{x}}{dt^2}\cdot\frac{d^2t}{d\tilde{x}^2} + \frac{d\tilde{x}}{dt} \cdot \frac{d^3\tilde{x}}{dt^3} \cdot  \frac{d^2t}{d\tilde{x}^2}  \nonumber \\
& &  +  \frac{d\tilde{x}}{dt}\cdot\frac{d^2\tilde{x}}{dt^2}\cdot \frac{d^3t}{d\tilde{x}^3} \cdot \frac{d\tilde{x}}{dt} \Bigg )
+ \left ( 3 \left ( \frac{d\tilde{x}}{dt} \right )^2 \cdot \frac{d^2\tilde{x}}{dt^2} \cdot \frac{d^3t}{d\tilde{x}^3} + \left ( \frac{d\tilde{x}}{dt} \right )^3\cdot \frac{d^4t}{d\tilde{x}^4} \cdot \frac{d\tilde{x}}{dt}  \right ) \nonumber \\
0 & = & \frac{d^4\tilde{x}}{dt^4}\cdot \frac{dt}{d\tilde{x}} + 4 \cdot \frac{d^3\tilde{x}}{dt^3} \cdot \frac{d\tilde{x}}{dt} \cdot \frac{d^2t}{d\tilde{x}^2}
+3 \cdot \frac{d^2\tilde{x}}{dt^2}\cdot\frac{d^2\tilde{x}}{dt^2} \cdot \frac{d^2t}{d\tilde{x}^2} \nonumber \\
& & + 6 \left ( \frac{d\tilde{x}}{dt} \right )^2 \cdot \frac{d^2\tilde{x}}{dt^2} \cdot \frac{d^3t}{d\tilde{x}^3} + \left ( \frac{d\tilde{x}}{dt} \right )^4 \cdot \frac{d^4t}{d\tilde{x}^4}
\end{eqnarray}

\vspace{.1in}

\noindent
If we keep taking the derivative with respect to $t$ for $k$ times and collect the similar terms after each differentiation as shown above, eventually we will arrive at an expression of the form

\begin{eqnarray}
\label {eq:x-tilde-t-k+2}
0 & = & \frac{d^{k+2}\tilde{x}}{dt^{k+2}} \cdot \frac{dt}{d\tilde{x}} + \left ( \text{terms involving } \frac{d^{k+1}\tilde{x}}{dt^{k+1}}, \frac{d^{k}\tilde{x}}{dt^{k}} , ... , \frac{d\tilde{x}}{dt} , \frac{dt}{d\tilde{x}},  \frac{d^2t}{d\tilde{x}^2}, ... , \frac{d^{k+1}t}{d\tilde{x}^{k+1}} \right ) \nonumber \\
& & + \left ( \frac{d\tilde{x}}{dt} \right )^{k+2} \cdot \frac{d^{k+2}t}{d\tilde{x}^{k+2}}
\end{eqnarray}

\noindent
From (\ref{eq:t_k+2}) we know that

\begin{eqnarray}
\label{eq:t-over-tilde-x}
\frac{dt}{d\tilde{x}} | _{t=0} & = & 1 \nonumber \\
\frac{d^2t}{d\tilde{x}^2} | _{t=0}  & = & 2b_2 \nonumber \\
\vdots & = & \vdots \\
\frac{d^{k+1}t}{d\tilde{x}^{k+1}} | _{t=0}  & = & (k+1)!b_{k+1} \nonumber \\
\frac{d^{k+2}t}{d\tilde{x}^{k+2}} | _{t=0} & = & (k+2)!b_{k+2} \nonumber
\end{eqnarray}

\vspace{.1in}

\noindent
Then we look at $\frac{d\tilde{x}}{dt} | _{t=0}  $, ... , $\frac{d^{k+1}\tilde{x}}{dt^{k+1}}| _{t=0} $, and $\frac{d^{k+2}\tilde{x}}{dt^{k+2}}| _{t=0} .$

\vspace{.1in}

\noindent
By definition (\ref {eq:tilde-x-defn}) and Equation (\ref {eq:u_m+k}),
\begin{eqnarray}
\tilde{x} & = & \left (\frac{at^m + a_{m+1}t^{m+1} +   a_{m+2}t^{m+2} + \cdots +  a_{m+k}t^{m+k} +  a_{m+k+1}t^{m+k+1} + O \left ( t^{m+k+2}\right )}{a} \right )^{\frac{1}{m}}  \nonumber \\
& = &  t \left ( 1 + \frac{a_{m+1}}{a} t + \frac{a_{m+2}}{a} t^2 +\cdots +  \frac{a_{m+k}}{a}t^{k} +  \frac{a_{m+k+1}}{a}t^{k+1} + O \left ( t^{k+2}\right ) \right )^{\frac{1}{m}}  \nonumber \\
& = & t \Bigg \{ 1+  \frac{1}{m} \left [ \frac{a_{m+1}}{a} t + \frac{a_{m+2}}{a} t^2 + \cdots +  \frac{a_{m+k}}{a}t^{k} +  \frac{a_{m+k+1}}{a}t^{k+1} + O \left ( t^{k+2}\right ) \right ] \nonumber \\
& & + \frac{1}{2}\cdot \frac{1}{m} \left ( \frac{1}{m} -1 \right ) \left [ \frac{a_{m+1}}{a} t + \frac{a_{m+2}}{a} t^2 + \cdots +  \frac{a_{m+k}}{a}t^{k} +  \frac{a_{m+k+1}}{a}t^{k+1} + O \left ( t^{k+2}\right ) \right] ^2 \nonumber \\
& & + \cdots \nonumber \\
& & + \frac{1}{(k+1)!} \cdot \frac{1}{m} \cdot \left ( \frac{1}{m} -1 \right ) \cdots \left ( \frac{1}{m} -k \right ) \Big [ \frac{a_{m+1}}{a} t + \frac{a_{m+2}}{a} t^2 + \cdots  +  \frac{a_{m+k}}{a}t^{k} \nonumber  \\
& & +  \frac{a_{m+k+1}}{a}t^{k+1} + O \left ( t^{k+2}\right ) \Big ] ^{k+1} + O \left ( t^{k+2}\right ) \Bigg \} \nonumber
\end{eqnarray}

\noindent
After collecting similar terms we can write
\begin{equation}
\label{eq:tilde-x-induction}
\tilde{x}= t + \lambda_2t^2 + \lambda_3t^3 + \cdots + \lambda_{k+1}t^{k+1} + \left ( \frac{a_{m+k+1}}{ma} + \lambda_{k+2} \right ) t^{k+2} + O \left ( t^{k+3} \right ),
\end{equation}
where

\begin{itemize}
	\item $\lambda_2$ is a constant involving $m, a,$ and $a_{m+1}$
	\item $\lambda_3$ is a constant involving $m, a, a_{m+1}$ and $a_{m+2}$
\end{itemize}
 	
	\vdots
	\begin{itemize}
	\item $\lambda_{k+1}$ is a constant involving $m, a, a_{m+1}, ..., a_{m+k-1}, a_{m+k}$
	\item $\lambda_{k+2}$ is a constant involving $m, a, a_{m+1}, ..., a_{m+k-1}, a_{m+k}$
\end{itemize}

\vspace{.1in}

\noindent
By the inductive hypothesis, $\lambda_2, \lambda_3, ... \lambda_{k+1}, \lambda_{k+2}$ are all constants that only depend on $m,n$, and the function $f$.

\vspace{.1in}

\noindent
From Equation (\ref {eq:tilde-x-induction}) we obtain

\begin{eqnarray}
\label{eq:tilde-x-over-t}
\frac{d\tilde{x}}{dt} | _{t=0}  & = & 1 \nonumber \\
\frac{d^2\tilde{x}}{dt^2} | _{t=0}  & = & 2\lambda_2 \nonumber \\
\vdots & = & \vdots \\
\frac{d^{k+1}\tilde{x}}{dt^{k+1}} | _{t=0}  & = & (k+1)!\lambda_{k+1} \nonumber \\
\frac{d^{k+2}\tilde{x}}{dt^{k+2}} | _{t=0}  & = & (k+2)! \left ( \frac{a_{m+k+1}}{ma} + \lambda_{k+2} \right ) \nonumber
\end{eqnarray}

\vspace{.1in}
\noindent
Now we evaluate (\ref {eq:x-tilde-t-k+2}) at $t=0$ and make use of (\ref {eq:t-over-tilde-x}) and (\ref {eq:tilde-x-over-t}):

$$ 0 = (k+2)! \left ( \frac{a_{m+k+1}}{ma} + \lambda_{k+2} \right ) \cdot 1 + \Big ( \text{terms involving } b_2, ..., b_{k+1}, \lambda_2, ..., \lambda_{k+1} \Big ) + 1 \cdot  (k+2)!b_{k+2} .$$

\noindent
Thus we obtain

\begin{equation}
\label{eq:b_k+2}
b_{k+2} = - \frac{a_{m+k+1}}{ma} + Q,
\end{equation}

\vspace{.1in}
\noindent
where $Q$ is a constant depending on $b_2, ..., b_{k+1}, \lambda_2, ..., \lambda_{k+1}, \lambda_{k+2}$, hence $Q$ is completely determined by $m,n$, and $f$.

\vspace{.1in}

\noindent
Next we will analyze $a_{m+k+1}$.  From (\ref {eq:t_k+2}) we have

\begin{eqnarray*}
t^{m-n} & = & \tilde{x}^{m-n} \left ( 1+ b_2 \tilde{x} +   \cdots + b_{k+1} \tilde{x}^{k} + b_{k+2} \tilde{x}^{k+1} + O \left (  \tilde{x}^{k+2} \right )\right ]^{m-n} \\
& = & \tilde{x}^{m-n} \Big \{ 1+ (m-n) \left [ b_2 \tilde{x}  + \cdots + b_{k+1} \tilde{x}^{k} + b_{k+2} \tilde{x}^{k+1} + O \left (  \tilde{x}^{k+2} \right ) \right ]   \\
& & +  \frac{(m-n)(m-n-1)}{2}  \left [ b_2 \tilde{x} +  \cdots + b_{k+1} \tilde{x}^{k} + b_{k+2} \tilde{x}^{k+1} + O \left (  \tilde{x}^{k+2} \right ) \right ]^2 + \cdots   \\
& &  + \frac{(m-n)(m-n-1) \cdots (m-n-k) }{(k+1)!}  \left [ b_2 \tilde{x} +  \cdots + b_{k+1} \tilde{x}^{k} + b_{k+2} \tilde{x}^{k+1} + O \left (  \tilde{x}^{k+2} \right ) \right ]^{k+1}   \\
& & + O \left ( \tilde{x} ^ {k+2}\right ) \Big \}
\end{eqnarray*}

\vspace{.1in}
\noindent
After collecting similar terms we can express $t^{m-n}$ as

\begin{eqnarray*}
t^{m-n} &  = &  \tilde{x}^{m-n} \Big \{ 1+ c_{1,m-n}\tilde{x} + c_{2,m-n} \tilde{x}^2 + \cdots + c_{k,m-n} \tilde{x}^k  \\
& & + \left [ (m-n)b_{k+2} + c_{k+1,m-n} \right ] \tilde{x}^{k+1} +  O \left (  \tilde{x}^{k+2} \right ) \Big \} ,
\end{eqnarray*}

\vspace{.1in}
\noindent
where $c_{1,m-n}$ is a constant depending on $m$ and $b_2$; $c_{2,m-n}$ is a constant depending on $m,b_2$ and $b_3$; ... ; $c_{k,m-n}$ is a constant depending on $m,b_2, ... , b_{k+1}$; $c_{k+1,m-n}$ is a constant depending on $m,b_2, ... , b_{k+1}$.

\vspace{.1in}

\noindent
By the inductive hypothesis, $c_{1,m-n}, c_{2,m-n}, ..., c_{k,m-n}$ and $c_{k+1,m-n}$ are all determined only by $m,n$, and $f$.  Thus we have

\begin{eqnarray}
\label{eq:t_m-n_induction}
t^{m-n}&  = & \tilde{x}^{m-n} + c_{1,m-n}\tilde{x}^{m-n+1} + c_{2,m-n} \tilde{x}^{m-n+2} + \cdots + c_{k,m-n} \tilde{x}^{m-n+k} \nonumber \\
& &  + \left [ (m-n)b_{k+2} + c_{k+1,m-n} \right ] \tilde{x}^{m-n+k+1} +  O \left (  \tilde{x}^{m-n+k+2} \right ) ,
\end{eqnarray}

\vspace{.1in}
\noindent
where $c_{1,m-n}, c_{2,m-n}, ..., c_{k,m-n}$ and $c_{k+1,m-n}$ are constants depending on $m,n$, and $f$.

\vspace{.1in}

\noindent
By the same type of analysis we obtain similar expressions for the other powers of $t$:

\vspace{.1in}

\begin{eqnarray}
\label{eq:t_m-n+1_induction}
t^{m-n+1} & = & \tilde{x}^{m-n+1} + c_{1,m-n+1}\tilde{x}^{m-n+2} + c_{2,m-n+1} \tilde{x}^{m-n+3} + \cdots + c_{k,m-n+1} \tilde{x}^{m-n+k+1}  \nonumber \\
& & + \left [ (m-n+1)b_{k+2} + c_{k+1,m-n+1} \right ] \tilde{x}^{m-n+k+2} +  O \left (  \tilde{x}^{m-n+k+3} \right ) ,
\end{eqnarray}

\vspace{.1in}
\noindent
where $c_{1,m-n+1}, c_{2,m-n+1}, ..., c_{k,m-n+1}$ and $c_{k+1,m-n+1}$ are constants depending on $m,n$, and $f$.

\vspace{.1in}

\begin{eqnarray}
\label{eq:t_m-n+2_induction}
t^{m-n+2} & = & \tilde{x}^{m-n+2} + c_{1,m-n+2}\tilde{x}^{m-n+3} + c_{2,m-n+2} \tilde{x}^{m-n+4} + \cdots + c_{k,m-n+2} \tilde{x}^{m-n+k+2}  \nonumber \\
& & + \left [ (m-n+2)b_{k+2} + c_{k+1,m-n+2} \right ] \tilde{x}^{m-n+k+3} +  O \left (  \tilde{x}^{m-n+k+4} \right ) ,
\end{eqnarray}
where $c_{1,m-n+2}, c_{2,m-n+2}, ..., c_{k,m-n+2}$ and $c_{k+1,m-n+2}$ are constants depending on $m,n$, and $f$.

\vdots

\begin{eqnarray}
\label{eq:t_m-n+k_induction}
t^{m-n+k} & = & \tilde{x}^{m-n+k} + c_{1,m-n+k}\tilde{x}^{m-n+k+1} + c_{2,m-n+k} \tilde{x}^{m-n+k+2} + \cdots + c_{k,m-n+k} \tilde{x}^{m-n+2k}  \nonumber \\
& & + \left [ (m-n+k)b_{k+2} + c_{k+1,m-n+k} \right ] \tilde{x}^{m-n+2k+1} +  O \left (  \tilde{x}^{m-n+2k+2} \right ) ,
\end{eqnarray}
where $c_{1,m-n+k}, c_{2,m-n+k}, ..., c_{k,m-n+k}$ and $c_{k+1,m-n+k}$ are constants depending on $m,n$, and $f$.

\vspace{.1in}

\begin{eqnarray}
\label{eq:t_m-n+k+1_induction}
t^{m-n+k+1} & = & \tilde{x}^{m-n+k+1} + c_{1,m-n+k+1}\tilde{x}^{m-n+k+2} \nonumber \\
& & + c_{2,m-n+k+1} \tilde{x}^{m-n+k+3} + \cdots  + c_{k,m-n+k+1} \tilde{x}^{m-n+2k+1}  \nonumber \\
& & + \left [ (m-n+k+1)b_{k+2} + c_{k+1,m-n+k+1} \right ] \tilde{x}^{m-n+2k+2} +  O \left (  \tilde{x}^{m-n+2k+3} \right ) ,
\end{eqnarray}
where $c_{1,m-n+k+1}, c_{2,m-n+k+1}, ..., c_{k,m-n+k+1}$ and $c_{k+1,m-n+k+1}$ are constants depending on $m,n$, and $f$.

\vspace{.1in}

\begin{equation}
\label{eq:t_m-n+k+2_induction}
t^{m-n+k+2} = O \left ( \tilde{x}^{m-n+k+2} \right )
\end{equation}

\vspace{.1in}

\noindent
From (\ref {eq:u_m+k}) we obtain

\begin{eqnarray*}
u^{(n)} & = & am(m-1) \cdots (m-n+1) t^{m-n} + a_{m+1}(m+1)m \cdots (m-n+2) t^{m-n+1} \\
& & + \cdots + a_{m+k}(m+k)(m+k-1) \cdots (m-n+k+1)t^{m-n+k} \\
& &  + a_{m+k+1}(m+k+1)(m+k) \cdots (m-n+k+2) t^{m-n+k+1} + O \left (  t^{m-n+k+2}\right )\\
\end{eqnarray*}

\noindent
Then by (\ref {eq:t_m-n_induction}) to (\ref {eq:t_m-n+k+2_induction}) we can write
\begin{eqnarray*}
u^{(n)} & = & am(m-1) \cdots (m-n+1) \Big \{ \tilde{x}^{m-n} + c_{1,m-n}\tilde{x}^{m-n+1} + c_{2,m-n} \tilde{x}^{m-n+2} + \cdots \\
& & + c_{k,m-n} \tilde{x}^{m-n+k}  + \left [ (m-n)b_{k+2} + c_{k+1,m-n} \right ] \tilde{x}^{m-n+k+1} +  O \left (  \tilde{x}^{m-n+k+2} \right ) \Big \} \\
& &  + a_{m+1}(m+1)m \cdots (m-n+2) \Big \{  \tilde{x}^{m-n+1} + c_{1,m-n+1}\tilde{x}^{m-n+2} + c_{2,m-n+1} \tilde{x}^{m-n+3} \\
& & + \cdots + c_{k,m-n+1} \tilde{x}^{m-n+k+1}  + \left [ (m-n+1)b_{k+2} + c_{k+1,m-n+1} \right ] \tilde{x}^{m-n+k+2} \\
& & +  O \left (  \tilde{x}^{m-n+k+3} \right ) \Big \}  + \cdots \\
& & + a_{m+k}(m+k)(m+k-1) \cdots (m+k-n+1) \Big \{  \tilde{x}^{m-n+k} + c_{1,m-n+k}\tilde{x}^{m-n+k+1} \\
& & + c_{2,m-n+k} \tilde{x}^{m-n+k+2} + \cdots + c_{k,m-n+k} \tilde{x}^{m-n+2k} +  \\
& &  \left [ (m-n+k)b_{k+2} + c_{k+1,m-n+k} \right ] \tilde{x}^{m-n+2k+1} +  O \left (  \tilde{x}^{m-n+2k+2} \right )  \Big \} \\
& & + a_{m+k+1}(m+k+1)(m+k) \cdots (m+k-n+2) \Big \{ \tilde{x}^{m-n+k+1} +\\
& & c_{1,m-n+k+1}\tilde{x}^{m-n+k+2}  + c_{2,m-n+k+1} \tilde{x}^{m-n+k+3} + \cdots  + c_{k,m-n+k+1} \tilde{x}^{m-n+2k+1}  \\
& & + \left [ (m-n+k+1)b_{k+2} + c_{k+1,m-n+k+1} \right ] \tilde{x}^{m-n+2k+2} +  O \left (  \tilde{x}^{m-n+2k+3} \right )  \Big \} \\
& & + O \left ( \tilde{x}^{m-n+k+2} \right ) \\
& = &  am(m-1) \cdots (m-n+1) \tilde{x}^{m-n} + C \left (m,a,a_{m+1},c_{1,m-n} \right ) \tilde{x}^{m-n+1} \\
& & + C \left (m,a,a_{m+1}, a_{m+2}, c_{1,m-n+1} , c_{2,m-n}\right ) \tilde{x}^{m-n+2} + \cdots    \\
& & +  C \left (m,a,a_{m+1}, ..., a_{m+k}, c_{k,m-n} , c_{k-1,m-n+1},..., c_{1,m-n+k-1}\right )  \tilde{x}^{m-n+k} \\
& & + \Big [ am(m-1)\cdots (m-n)b_{k+2} + a_{m+k+1}(m+k+1)(m+k) \cdots (m+k-n+2) \\
& & + C \left (m,a,a_{m+1}, ..., a_{m+k}, c_{k+1,m-n} , c_{k,m-n+1},..., c_{1,m-n+k}\right ) \Big ] \tilde{x}^{m-n+k+1}  \\
& & + O \left ( \tilde{x}^{m-n+k+2} \right ).
\end{eqnarray*}

\vspace{.1in}

\noindent
Here $C \left (m,a,a_{m+1},c_{1,m-n} \right )$ is a constant depending on  $m,a,a_{m+1}$ and $c_{1,m-n}$; we denote it as $p_{m-n+1}$ to simplify notations. Since $a,a_{m+1}$ and $c_{1,m-n}$ only depend on $m,n$, and $f$, we know that $p_{m-n+1}$ only depends on $m,n$, and $f$.

\vspace{.1in}

\noindent
Similarly, the other constants $C \left (m,a,a_{m+1}, a_{m+2}, c_{1,m-n+1} , c_{2,m-n}\right )  ,$
......, and
\newline
$C \left (m,a,a_{m+1}, ..., a_{m+k}, c_{k+1,m-n} , c_{k,m-n+1},..., c_{1,m-n+k}\right )$ all depend only on $m,n$, and $f$, and can be denoted simply as $p_{m-n+2}$,..., $p_{m-n+k}$, and $p_{m-n+k+1}$.   Thus we can rewrite the above equation as

\begin{eqnarray}
\label{eq:u^(n)_induction}
 u^{(n)} & = & am(m-1) \cdots (m-n+1) \tilde{x}^{m-n} + p_{m-n+1} \tilde{x}^{m-n+1} + \nonumber \\
& &  p_{m-n+2} \tilde{x}^{m-n+2} + \cdots   + p_{m-n+k}   \tilde{x}^{m-n+k} + \nonumber \\
& & \Big [ am(m-1)\cdots (m-n)b_{k+2} + a_{m+k+1}(m+k+1)(m+k) \cdots (m+k-n+2) \nonumber \\
& & + p_{m-n+k+1}  \Big ] \tilde{x}^{m-n+k+1}  + O \left ( \tilde{x}^{m-n+k+2} \right )
\end{eqnarray}

\noindent
Now because of $u^{(n)}=f(u)$ and definition (\ref {eq:tilde-x-defn}) , we have

\begin{eqnarray}
f(u) & = & am(m-1) \cdots (m-n+1) \left (  \frac{u}{a}\right )^{\frac{m-n}{m}} + p_{m-n+1}  \left (  \frac{u}{a}\right )^{\frac{m-n+1}{m}} + \nonumber \\
& &  p_{m-n+2} \left (  \frac{u}{a}\right )^{\frac{m-n+2}{m}} + \cdots   + p_{m-n+k}  \left (  \frac{u}{a}\right )^{\frac{m-n+k}{m}} + \nonumber \\
& &  \Big [ am(m-1)\cdots (m-n)b_{k+2} + a_{m+k+1}(m+k+1)(m+k) \cdots (m+k-n+2) \nonumber \\
& & + p_{m-n+k+1}   \Big ] \left (  \frac{u}{a}\right )^{\frac{m-n+k+1}{m}}  + O \left ( u^{\frac{m-n+k+2}{m}} \right ) \nonumber
\end{eqnarray}

\vspace{.1in}
\noindent
Due to (\ref {eq:b_k+2}), we can rewrite the above equation as

\begin{eqnarray}
\label{eq:f(u)_induction}
f(u) & = &  am(m-1) \cdots (m-n+1) \left (  \frac{u}{a}\right )^{\frac{m-n}{m}} + p_{m-n+1}  \left (  \frac{u}{a}\right )^{\frac{m-n+1}{m}} + \nonumber \\
& &  p_{m-n+2} \left (  \frac{u}{a}\right )^{\frac{m-n+2}{m}} + \cdots   + p_{m-n+k}  \left (  \frac{u}{a}\right )^{\frac{m-n+k}{m}} + \nonumber \\
& & + \Bigg \{ \Big [ (m+k+1)(m+k) \cdots (m+k-n+2)- (m-1)(m-2) \cdots (m-n)\Big ] a_{m+k+1} \nonumber \\
& & + am(m-1)\cdots (m-n)Q  + p_{m-n+k+1}  \Bigg \} \left (  \frac{u}{a}\right )^{\frac{m-n+k+1}{m}}   + O \left ( u^{\frac{m-n+k+2}{m}} \right )
\end{eqnarray}

\noindent
From (\ref {eq:f(u)_induction}) we get

\begin{eqnarray}
\label{eq:a_m+k+1_induction}
& & \Big [ (m+k+1)(m+k) \cdots (m+k-n+2)- (m-1)(m-2) \cdots (m-n)\Big ] a_{m+k+1} \nonumber \\
& & + am(m-1)\cdots (m-n)Q  + p_{m-n+k+1} \nonumber \\
& =  & \lim _{u \to 0} \frac{f(u)-am \cdots (m-n+1) \left (  \frac{u}{a}\right )^{\frac{m-n}{m}} - p_{m-n+1}  \left (  \frac{u}{a}\right )^{\frac{m-n+1}{m}}  - \cdots   - p_{m-n+k}  \left (  \frac{u}{a}\right )^{\frac{m-n+k}{m}}  }{\left (  \frac{u}{a}\right )^{\frac{m-n+k+1}{m}}} \nonumber
\end{eqnarray}

\vspace{.1in}

\noindent
Note that $(m+k+1)(m+k) \cdots (m+k-n+2)- (m-1)(m-2) \cdots (m-n) \neq 0$, then since the constants $Q, a, p_{m-n+1},..., p_{m-n+k+1}$ all depend only on $m,n$, and $f$, we know that $a_{m+k+1}$ only depends on $m,n$, and $f$.  Consequently, $b_{k+2}$ also only depends on $m,n$, and $f$ because of (\ref {eq:b_k+2}).

\vspace{.1in}

\noindent
Therefore, by mathematical induction, all derivatives of $f$ at 0 are determined completely by $m,n$, and $f$.  This completes the proof of Lemma \ref{lem:step1_Taylor_coefficient}.

\vspace{.15in}

\section {The Proofs of Theorems \ref{thm:uniqueness-ODE}, \ref {thm:u1_u2} and \ref{thm:f_ho"lder}}
\label{sec:main_pf}

\noindent
\textbf{Theorem \ref{thm:uniqueness-ODE} }

\vspace{.05in}
\pf
\noindent
Suppose there are two solutions $u(t)$ and $v(t)$, both satisfy Equations (\ref{eq:main}) and  (\ref{eq:initial}).  By Lemma \ref{lem:step1_Taylor_coefficient}, at $t=0$, $u$ and $v$ have the same derivative of any order.  Let $w=u-v$, then $$w^{(k)}(0)=0 \hspace{.2in}  \text {for any integer} \hspace{.1in} k \geq 0.$$  In order to apply Theorem \ref{thm:Pan-Wang_flatness} we need to show that $w$ satisfies Condition (\ref{eq: 2nd_derivative_condition}).

$$
	w^{(n)}(t) \,\, = \,\, u^{(n)}(t)-v^{(n)}(t) \,\, = \,\,  f(u(t))-f(v(t))\\
$$

\noindent
Without loss of generality we assume $a>0$.  By Equation (\ref {eq:for-use-in-main-thm}), we can write
\begin{equation}
\label {eq:for-f-differentiability}
 f(u) = \left [ a^{\frac{n}{m}}m(m-1) \cdots (m-n+1) u^{\frac{m-n}{m}}+\alpha (u) \right ]
\end{equation}
and

$$f(v) = \left [ a^{\frac{n}{m}}m(m-1) \cdots (m-n+1) v^{\frac{m-n}{m}}+\alpha (v) \right ],$$
where $\alpha$ is a function with the order $$\alpha(s)=O \left ( s^{\frac{m-n+1}{m}}\right ).$$

\noindent
So we can write
\begin{equation}
\label{eq:w^(n)(t)}
w^{(n)}(t)=a^{\frac{n}{m}}m(m-1) \cdots (m-n+1) \left ( u^{\frac{m-n}{m}} - v^{\frac{m-n}{m}}  \right ) + \left ( \alpha (u) - \alpha (v) \right )
\end{equation}

\vspace{.1in}

\noindent
If $m=n$, then $$a^{\frac{n}{m}}m(m-1) \cdots (m-n+1) \left ( u^{\frac{m-n}{m}} - v^{\frac{m-n}{m}}  \right )=0$$

\vspace{.1in}

\noindent
If $m>n$, by the Mean Value Theorem,
\begin{equation}
\label{eq:u-and-v-exponent}
  \left | u^{\frac{m-n}{m}} - v^{\frac{m-n}{m}}   \right | \leq \frac{m-n}{m}\zeta^{-\frac{n}{m}} |u-v|
	\end{equation}
	where $\zeta(t)$ is between $u(t)$ and $v(t)$.  Since $u(t)=at^m+O \left ( t^{m+1} \right ) $ and $v(t)=at^m+O \left ( t^{m+1} \right )$, we know that $\zeta(t)=at^m+O \left ( t^{m+1} \right )$, which implies
$$\zeta ^{\frac{n}{m}} \,\, = \,\, a^{\frac{n}{m}} t^n \left (  1+O(t) \right ) \,\, \geq \,\, Ct^n$$ for some constant $C>0$ when $t$ is sufficiently small.
Thus $$\zeta^{-\frac{n}{m}} |u-v| \leq C^{-1} \frac{|u-v|}{t^n}$$ and by equation (\ref{eq:u-and-v-exponent}) we know that
\begin{equation}
\label{eq:first-part-estimate}
a^{\frac{n}{m}}m(m-1) \cdots (m-n+1)\left | u^{\frac{m-n}{m}} - v^{\frac{m-n}{m}}  \right | \leq C \frac{|u-v|}{t^n}
\end{equation}
for another constant $C>0$.

\vspace{.1in}
\noindent
Next, we estimate $\left | \alpha (u) - \alpha (v) \right | $.

\vspace{.1in}
\noindent
From Equation (\ref {eq:main}), we know that $f$ is differentiable with respect to $t$, since $u^{(n)}(t)$ is differentiable with respect to $t$.  Condition (\ref {eq:initial}) shows that $\frac{du}{dt} (t) \neq 0$ when $t>0$ is sufficiently small.  Then by the Inverse Function Theorem, $t$ is differentiable with respect to $u$.  Thus, when $u$ is small and positive, $f$ is differentiable with respect to $u$ and $$\frac{df}{du}=\frac{df}{dt}\cdot \frac{dt}{du}.$$

\vspace{.1in}

\noindent
Then by Equation (\ref {eq:for-f-differentiability}), since $f$ is differentiable on a small interval $(0,\delta)$, $\alpha$ is also differentiable on a small interval $(0,\delta)$.  By the Mean Value Theorem

$$ \alpha (u) - \alpha (v) = \alpha '(\eta) (u-v) $$ where $\eta(t)$ is between $u(t)$ and $v(t)$.  Because $u(t)=at^m+O \left ( t^{m+1} \right ) $ and $v(t)=at^m+O \left ( t^{m+1} \right )$, we know that $\eta(t)=at^m+O \left ( t^{m+1} \right ) \geq Ct^m$ when $t$ is small, thus $$\eta ^{-\frac{n-1}{m}}=O \left (t^{-(n-1)} \right ).$$

\noindent
From  $\alpha (s) = O \left ( s^{\frac{m-n+1}{m}}\right ) $ we get $\alpha ' (s) = O \left ( s^{-\frac{n-1}{m}}\right ) $.  Therefore $$\alpha '(\eta) \,\, = \,\, O \left ( \eta^{-\frac{n-1}{m}}\right ) \,\, = \,\, O \left (t^{-(n-1)} \right ) .$$
Thus for some $C>0$

\begin{eqnarray}
\label{eq:2nd-part-estimate}
\left | \alpha (u) - \alpha (v) \right | & \leq & C t^{-(n-1)} |u-v| \nonumber \\
&  \leq &  C t^{-n} |u-v|  \hspace{.2in} \text{ since} \hspace{.1in} 0<t<1
\end{eqnarray}

\noindent
Combining Equations (\ref{eq:w^(n)(t)}), (\ref{eq:first-part-estimate}), and (\ref{eq:2nd-part-estimate}), we conclude $$|w^{(n)}(t)| \,\, \leq \,\, C \frac{|u(t)-v(t)|}{t^n} \,\, = \,\,  C \frac{|w(t)|}{t^n}.$$

\vspace{.1in}

\noindent
Finally, extend the domain of $w(t)$ to $[-1,1]$ by defining $w(t)=w(-t)$ when $-1 \leq t <0$.  Then $w \in C^{\infty} ([-1,1])$ and it satisfies Condition (\ref{eq: 2nd_derivative_condition}).  By Theorem \ref{thm:Pan-Wang_flatness}, $w \equiv 0$, which means $u \equiv v$.

\noindent
This completes the proof of Theorem \ref{thm:uniqueness-ODE}.

\stop

\vspace{.1in}

\noindent
\textbf{Theorem \ref{thm:u1_u2}}

\vspace{.05in}
\pf  Without loss of generality we assume $a>0$.  We apply the same analysis as in the proof of Lemma \ref {lem:step1_Taylor_coefficient} to $u_1$ and $u_2$, respectively.  Similar to (\ref  {eq:a-determined}) we have
$$a^{\frac{n}{m}} \,\, = \,\, \lim _{u_1 \to 0} \frac{f(u_1)}{m(m-1) \cdots (m-n+1) u_1^{\frac{m-n}{m}}} \,\, = \,\, \lim _{s \to 0} \frac{f(s)}{m(m-1) \cdots (m-n+1) s^{\frac{m-n}{m}}} $$
and
$$b^{\frac{n}{l}} \,\, = \,\, \lim _{u_2 \to 0} \frac{f(u_2)}{l(l-1) \cdots (l-n+1) u_2^{\frac{l-n}{l}}} \,\, = \,\, \lim _{s \to 0} \frac{f(s)}{l(l-1) \cdots (l-n+1) s^{\frac{l-n}{l}}}.$$

\vspace{.1in}
\noindent
Suppose $m \neq l$, without loss of generality we assume $m<l$. Dividing the two equations, we get

\begin{eqnarray*}
a^{ \frac{n}{m}}b^{ -\frac{n}{l}} & = & \frac{l(l-1)\cdots (l-n+1)}{m(m-1)\cdots (m-n+1)} \lim _{s \to 0} s^{\frac{l-n}{l}-\frac{m-n}{m}} \\
& = & \frac{l(l-1)\cdots (l-n+1)}{m(m-1)\cdots (m-n+1)} \lim _{s \to 0} s^{\frac{n}{m}-\frac{n}{l}}.
\end{eqnarray*}

\noindent
Since $m < l$, $\displaystyle \lim _{s \to 0} s^{\frac{n}{m}-\frac{n}{l}}=0$.  However, $a^{ \frac{n}{m}}b^{ -\frac{n}{l}} \neq 0$.  This is a contradiction.

\vspace{.1in}

\noindent
Therefore $m=l$, consequently $a=b$.  Then by Theorem \ref {thm:uniqueness-ODE} we know that $u_1 \equiv u_2$ for small $t$.

\stop

\vspace{.1in}

\noindent
\textbf{Theorem \ref{thm:f_ho"lder}:}

\vspace{.05in}
\pf: The proof of Lemma \ref{lem:step1_Taylor_coefficient} shows that near 0, $t$ is a function of $u$, therefore $u^{(n)}(t)$ can be expressed as a function $f$ of $u$.  Thus Equation (\ref{eq:main}) holds when $t>0$ is small.  From Condition (\ref {eq:initial}) we define $f(0)=0$.

\vspace{.05in}
\noindent
By the first two equations in (\ref {eq:for-use-in-main-thm}) and the discussions in Appendix \ref{sec:appendix}, we know that there is a function $h$ which is $C^1$ on the closed interval $[0, \epsilon]$ for some $\epsilon >0$, such that

$$f(u)  = am(m-1)\cdots(m-n+1)\tilde{x}^{m-n} + h \left ( \tilde{x} \right ) \tilde{x}^{m-n}.$$
By definition (\ref {eq:tilde-x-defn}) we have
\begin{equation}
\label{eq:f_Ho'lder-expression}
f(u)  =   am(m-1) \cdots(m-n+1)\left ( \frac{u}{a} \right ) ^ {\frac{m-n}{m}}+ h \left ( \left (\frac{u}{a} \right )^{\frac{1}{m}}\right ) \left ( \frac{u}{a} \right )^\frac{m-n}{m}
\end{equation}

\noindent
Since $0 \leq \frac{m-n}{m} <1$ and $0<\frac{1}{m}<1$ , it is well known that $u^{\frac{m-n}{m}}$ and $u^{\frac{1}{m}}$ are H\"older continuous on the closed interval $[0,1]$ with H\"older coefficients $\frac{m-n}{m}$ and $\frac{1}{m}$, respectively.  This implies that the first term in (\ref {eq:f_Ho'lder-expression}) is H\"older continuous on $[0,1]$.

\vspace{.05in}

\noindent
Since $h$ is $C^1$ on $[0, \epsilon]$, it is also H\"older continuous on $[0,\epsilon]$.  Then since the composition of two H\"older continuous functions is H\"older continuous, we know that $h \left ( \left (\frac{u}{a} \right )^{\frac{1}{m}}\right )$ is H\"older continuous with respect to $u$ on a closed interval $[0,\delta]$ with $\delta>0$.  Next, because the product of two H\"older continuous functions is also H\"older continuous, we know that $h \left ( \left (\frac{u}{a} \right )^{\frac{1}{m}}\right ) \cdot \left ( \frac{u}{a} \right )^\frac{m-n}{m} $ is H\"older continuous.  Thus the second term in (\ref {eq:f_Ho'lder-expression}) is H\"older continuous on $[0,\delta]$.

\vspace{.05in}

\noindent
Therefore, $f$ is H\"older continuous on $[0,\delta]$ and the theorem is proved.

\stop

\vspace{.05in}

\appendix \section {Differentiation of the Taylor Expansion}
\label{sec:appendix}

We will discuss the regularity of the remainder term in the Taylor expansion of a function that is used in the proof of Theorem \ref{thm:f_ho"lder} and the differentiation of the Taylor expansion that is frequently used in the proof of Lemma \ref{lem:step1_Taylor_coefficient}.

\vspace{.1in}

\noindent
In general, suppose a function $g(x) \in C^{k+1} \left ([a,b] \right )$, by the Taylor Theorem we can write

\begin{equation}
\label{eq:Taylor-expansion}
 g(x)=g(a)+g'(a)(x-a)+\frac{g''(a)}{2!}(x-a)^2+ \cdots + \frac{g^{(k)}(a)}{k!}{(x-a)^k} + h(x)(x-a)^k
\end{equation}
 where $\displaystyle \lim _{x \to a} h(x)=0$.  An explicit expression for $h(x)$ is
\begin{eqnarray}
\label{eq:Taylor_remainder}
h(x)=\frac{g^{(k+1)}(\xi)}{(k+1)!}(x-a)
\end{eqnarray}
where $a<\xi<x$.

\vspace{.1in}
\noindent
From (\ref {eq:Taylor-expansion}) we know that $h(x)$ is $C^1$ on $(a,b]$.  Next we show that it is actually $C^1$ up to the boundary, on $[a,b]$.
\vspace{.1in}

\noindent
Taking the derivative on both sides of Equation ( \ref {eq:Taylor-expansion} ), we get

\begin{equation}
\label{eq:g'(x)}
g'(x)=g'(a)+g''(a)(x-a)+ \cdots + \frac{g^{(k)}(a)}{(k-1)!}{(x-a)^{k-1}} + h'(x)(x-a)^k+k h(x) (x-a)^{k-1}.
\end{equation}

\vspace{.1in}

\noindent
Define $h(a)=0$, so $h$ is continuous on $[a,b]$.  Denote $$P(x)=g(a)+g'(a)(x-a)+\frac{g''(a)}{2!}(x-a)^2+ \cdots + \frac{g^{(k)}(a)}{k!}{(x-a)^k},$$ then $$h(x)=\frac{g(x)-P(x)}{(x-a)^k}.$$By the definition of limits,

\begin{eqnarray}
\label{eq:h-derivative-0}
h'(a) & = & \lim _{x \to a} \frac{h(x)-h(a)}{x-a} \nonumber \\
& = & \lim _{x \to a}\frac{g(x)-P(x)}{(x-a)^{k+1}} \nonumber \\
\vdots  \nonumber \\
& = & \frac{g^{(k+1)}(a)-P^{(k+1)}(a)}{(k+1)!} \hspace{.3in} \text{by applying the L'Hospital Rule} \,\, k+1 \,\, \text{ times} \nonumber \\
& = &  \frac{g^{(k+1)}(a)}{(k+1)!},
\end{eqnarray}
where we have used the fact that $P^{(k+1)}(a)=0$.

\vspace{.05in}
\noindent
When $x>a$,

\begin{eqnarray*}
h'(x) & = & \frac{d}{dx} \left (  \frac{g(x)-P(x)}{(x-a)^k} \right ) \\
& = & \frac{\left ( g'(x) - P'(x) \right )(x-a)^k - (g(x)-P(x))k(x-a)^{k-1} }{(x-a)^{2k}} \\
& = & \frac{g'(x) - P'(x)}{(x-a)^k} - \frac{k \left (  g(x)-P(x) \right )}{(x-a)^{k+1}}
\end{eqnarray*}

\noindent
By repeatedly applying the L'Hospital Rule, we know that

$$\lim _{x \to a} \frac{g'(x) - P'(x)}{(x-a)^k} \,\, = \,\, \frac{g^{(k+1)}(a)- P^{(k+1)}(a) }{k!} \,\, = \,\, \frac{g^{(k+1)}(a)}{k!}$$
and
$$\lim _{x \to a}  \frac{k \left (  g(x)-P(x) \right )}{(x-a)^{k+1}} \,\, = \,\, \frac{k \left (g^{(k+1)}(a)- P^{(k+1)} (a) \right )} {(k+1)!} \,\, = \,\, \frac{kg^{(k+1)}(a)}{(k+1)!}.$$

\noindent
Therefore

\begin{eqnarray}
\label{eq:h'(x)_limit}
\lim _{x \to a} h'(x) & = & \frac{g^{(k+1)}(a)}{k!} - \frac{kg^{(k+1)}(a)}{(k+1)!} \nonumber \\
& = & \frac{g^{(k+1)}(a)}{(k+1)!}
\end{eqnarray}

\vspace{.1in}
\noindent
Equations (\ref {eq:h-derivative-0}) and (\ref {eq:h'(x)_limit}) show that $h(x)$ is $C^1$ on the closed interval $[a,b]$.

\vspace{.15in}

\noindent
Furthermore, we know that for any $x \in [a,b]$, $ \left | h'(x) \right | \leq C$ for some constant $C_1$, thus $$ \left |  h'(x)(x-a)^k \right | \leq C_1 \left |  x-a \right |^k.$$ Since $g(x) \in C^{k+1} \left ([a,b] \right )$, from (\ref{eq:Taylor_remainder}) we know that $ \left | h(x) \right  | \leq C_2 \left |x-a \right |$ for some constant $C_2$, thus $$ \left |k h(x) (x-a)^{k-1} \right | \leq kC_2 \left |  x-a \right |^k.$$  Therefore, (\ref {eq:g'(x)}) can be denoted as

\begin{equation}
\label{eq:g'_O}
g'(x)=g'(a)+g''(a)(x-a)+ \cdots + \frac{g^{(k)}(a)}{(k-1)!}{(x-a)^{k-1}} + O \left ( x-a \right )^k
\end{equation}

\vspace{.1in}
\noindent
Since the first, second, ... , and $(k-1)$-th derivatives of $g'(x)$ at $a$ are $g''(a), g^{(3)}(a),$ ... , and $g^{(k)}(a)$, respectively, Equation (\ref{eq:g'_O}) is the Taylor expansion of $g'(x)$ at $a$ to order $k-1$.

\vspace{.05in}

\noindent
Usually we denote (\ref {eq:Taylor-expansion}) as

\begin{equation}
\label{eq:g_O}
g(x)=g(a)+g'(a)(x-a)+\frac{g''(a)}{2!}(x-a)^2+ \cdots + \frac{g^{(k)}(a)}{k!}{(x-a)^k} + O \left ( (x-a)^{k+1} \right ).
\end{equation}

\vspace{.05in}

\noindent
 This shows that we can formally differentiate (\ref{eq:g_O}) to get (\ref{eq:g'_O}).

\vspace{.2in}

\noindent
\textbf{Acknowledgment}: The research of Yu Yan was partially supported by the Scholar in Residence program at Indiana University-Purdue University Fort Wayne.

\vspace{.2in}

\bibliographystyle{plain}
\bibliography{thesis}

\end{document}